\documentclass[12pt]{article}
\usepackage[utf8]{inputenc}
\usepackage[all]{xy}
\usepackage{amsmath,amssymb}
\usepackage{enumerate}  

\def\cal#1{\mathcal{#1}}

\usepackage[all]{xy}

\usepackage[margin=2cm]{geometry}

\newtheorem{Th}{Theorem}

\newtheorem{Co}{Conjecture}

\author{Vladimir Blinovsky}

 \date{Instituto de Matematica e Statistica, USP,\\ Rua do Matao 1010, 05508- 090, Sao Paulo, Brazil;\\
Instituto de Ci\^encia e Tecnologia -- Unifesp, Avenida Cesare Mansueto Giulio Lattes, \\ 1201,  12247-014 
S\~ao Jos\'e dos Campos, SP, Brazil;\\
Institute for Information Transmission Problems 
B. Karetnyi 19, Moscow, Russia
\\ vblinovs@yandex.ru}

\begin{document}

\title{Tight Asymptotic of Probability of Singularity of $n\times n$ Random Matrix with Uniform Distributed $\pm 1$ Entries }

\maketitle
\bigskip

\begin{center}
{\bf Abstract}
\bigskip

We prove the conjecture: probability that $P_n$ of Bernulli $\pm 1$ square matrix is singular has tight asymptotic $4{n\choose 2}2^{-n}$. 
We also prove precise asymptotic 
$P_n -4{n\choose 2} 2^{-n} \sim 16 {n\choose 4}\left(\frac{3}{8}\right)^n$.
\end{center}
\bigskip

There is a number of works devoted to determination of tight asymptotic of probability $P_n$ that random square matrix with independent uniformly distributed $\pm 1$ entries is singular.

In this article~\cite{7} was stated the general
\begin{Co}
\label{co1}
The following asymptotic equality is valid
\begin{equation}
\label{e1}
P_n \sim 4{n\choose 2} 2^{-n}.
\end{equation}
\end{Co}
In this article we prove this conjecture.
\begin{Th}
\label{th1}
Asymptotic~(\ref{e1}) for $P_n$ is true.
\end{Th}
Using the same arguments as in the proof of  Conjecture~\ref{co1}, we prove
\begin{Co}
\label{co2}
The following asymptotic equality is valid
\begin{equation}
\label{e22}
P_n - 4{n\choose 2} 2^{-n}\sim16{n\choose 4} \left(\frac{3}{8}\right)^n .
\end{equation}
\end{Co}

\bigskip
{\bf History of the problem}.
\bigskip

Denote ${Q\choose i}$ the set of $i$- element subsets of finite set $Q$ and
$$
 [n]=\{1,2,\ldots, n\},\ [a,b]=\{a,\ldots ,b\},\ a\leq b\in \{1,2,\ldots\}.
 $$
Denote $P_n$ the probability that $n\times n$  random matrix with uniform distributed entries $\pm 1$ is singular.
Obvious lower bound for the value $P_n$ is the probability that two  or four rows or that two  or four 
columns of the matrix are linear dependent:
\begin{equation}
\label{l1}
P_n\geq 4{n\choose 2} 2^{-n}+16{n\choose 4} \left(\frac{3}{8}\right)^n -  12{n\choose 2}^2 4^{-n}  -2^{12} n^5 \left(\frac{3}{16}\right)^n .\end{equation}
Indeed following~\cite{cc1},  denote $\Gamma_1 ={[n]\choose 2}\times \{\pm 1\}\times \{L,R\},\ \Gamma_2 ={[n]\choose 4}\times\{\pm 1\}^3 \times \{L,R\}.$
So that $\alpha\in\Gamma_1$ specifies a set
of two distinct indices along with sign and direction bits. For given matrix $A_n$ with rows  $\{a_i , i\in [n]\}$ and columns $\{\bar{a}_i,\ i\in [n]\}$  event
$A_\alpha$ corresponds to the occurrence of a null vector of the form $\alpha$. For
example, $\alpha = (\{3, 6\}, -, R)\in\Gamma_1$ and $A_\alpha$ is the event that  $\bar{a}_3-\bar{a}_6 =0$ and if $\alpha =(\{3,5,7,8\}, -++,L)\in\Gamma_2$, then the event $B_\alpha$ is $a_3-a_5+a_7+a_8=0$. Denote $A=\bigcup_{\alpha\in \Gamma_1} A_\alpha,\ B=\bigcup_{\alpha\in \Gamma_2}B_\alpha$ 
and $W_1=\sum_{\alpha\in \Gamma_!}\hbox{Id}_{A_\alpha},\ W_2 =\sum_{\alpha\in \Gamma_2}\hbox{Id}_{B_\alpha}$, where   $\hbox{Id}_{A_\alpha}\ (\hbox{Id}_{B_\alpha})$ is indicator of the event $A_\alpha,\ \alpha\in \Gamma_1,\  (B_\alpha,\ \alpha\in \Gamma_2)$. We also denote $\Gamma_i (L) (\Gamma_i (R))$ when specifying the direction bits in $\Gamma_i$.

The inclusion-exclusion formula states that
\begin{equation}
\label{kkk100}
P(A)= \sum_{i=1}^{|\Gamma_1|} (-1)^{i+1} E{W_1\choose i},\ P(B)= \sum_{i=1}^{|\Gamma_2|} (-1)^{i+1} E{W_2\choose i},
\end{equation}
where ${W_m\choose i}=\frac{1}{i!}\prod_{j=1}^{i}(W_m-j+1),\ m=1,2$. 

The Bonferroni's inequalities states, that
\begin{eqnarray*}&&
P(A)\geq \sum_{i=1}^{|I|}(-1)^{i+1}E{W_1\choose i},\ I\subset \Gamma_1,\ |I|=2m,\\
&&
P(A)\leq \sum_{i=1}^{|I|}(-1)^{i+1}E{W_1\choose i},\ I\subset \Gamma_1,\ |I|=2m+1.
\end{eqnarray*}
The same inequalities we will use to estimate $P(B)$. 

We demonstrate the proof of Bonferroni's inequalities. If $\omega\in A$ is included in $r$ sets $A_\alpha$, then it counted  $\gamma=\sum_{\ell=1}^{|I|}(-1)^{\ell+1}{r\choose \ell}$ times in the rhs
of inequalities~(\ref{kkk100}). Here ${r\choose k}=0,\ k>r$. In the case when $k$ is odd, then $\gamma \geq 0$ and $\gamma \leq 0$ if $k$ is even. The remark that if $\omega\not\in A$,
then this event not evaluate in the rhs of inequalities complete the proof. 

It follows
\begin{equation}
\label{w1}
P(A)\geq E(W_1)-E{W_1\choose 2},\  P(B)\geq E(W_2)-E{W_2\choose 2}.
\end{equation}
 
We have
\begin{equation}
\label{v1}
E(W_1)=|\Gamma_1|2^{-n}=4{n\choose 2}2^{-n},\ E(W_2)=|\Gamma_2|\left(\frac{3}{8}\right)^n =16{n\choose 4}\left(\frac{3}{8}\right)^n . 
\end{equation}
In the first equality is counted the average number of pairs of linear dependent rows and pairs of linear dependent columns
in $A_n$. In second equality  is counted the average number quaternaries of linear dependent rows and  quaternaries  of linear dependent columns
in $A_n$. 
\begin{eqnarray}
&& \label{kkk101}
E{W_1\choose 2} =\frac{1}{2}(E(W_1^2)-E(W_1))  =\sum_{\alpha\neq \beta \in\Gamma_1 (L)} P(A_\alpha  A_\beta)+\sum_{\alpha\in\Gamma_1 (L),\beta\in\Gamma_1(R)} P(A_\alpha A_{\beta}) \\&& \nonumber
=4{n\choose 2}^2 2^{-2n}-4{n\choose 2}2^{-2n}+4{n\choose 2}^2 2^{-2n} 2 =  \left(12{n\choose 2}^2 -4{n\choose 2}\right)4^{-n}.
\end{eqnarray}
First summand in the second line  obtained by considering all pairs on one side $L$ or $R$ and, in second summand  we take away all pairs that share two rows along with $\pm$ combinations. The third  term has factor
of $2$, since for $\alpha,\ \beta$  on opposite sides, we have
$P(A_\alpha A_\beta) = 2P(A_\alpha )P(A_\beta )$. 

Next we have 
\begin{eqnarray}
\label{v2}
&& E{W_2 \choose 2}=\frac{1}{2}(E(W_2^2)-E(W_2))=\sum_{\alpha\neq \beta \in\Gamma_2 (L)} P(B_\alpha B_\beta)+\sum_{\alpha\in\Gamma_2 (L),\beta\in\Gamma_2(R)} P(B_\alpha B_{\beta})\\ \nonumber &&
\leq  2^{6}{n\choose 4}^2 \left(\frac{3}{8}\right)^{2n}+2^8 n^5\left(\frac{3}{16}\right)^n+ 2^6 \left(\frac{8}{3}\right)^8{n\choose 4}^2 \left(\frac{3}{8}\right)^{2n}<2^{16} n^5 \left(\frac{3}{16}\right)^{n}. \end{eqnarray}
Rhs. of  first inequality is evaluation of all pairs of quaternaries on one side $L$ or $R$ or on both sides. 
First plus second summands in the r.h.s. of first inequality  is the upper bound for the first sum in the lhs of the first inequality - this expression is just the square of $\sum_{\alpha\in\Gamma_2 (L)} P(B_\alpha)$ plus the upper bound for the evaluation of intersections of quaternaries on one side.  Third summand in the rhs of the first inequality  is the upper bound for the second sum in the lhs of the first inequality - this expression arises from the fact that besides the intersection of  quaternaries rows and columns -  submatrix of size $4\times 4$  another  elements of these quaternaries
are independent and hence for $\alpha\in\Gamma_2 (L),\ \beta\in\Gamma_2 (R)$ we have
$$
P(B_\alpha B_\beta) \leq \left(2^3\left( \frac{3}{8}\right)^{n-4}\right)^2.
$$

Then  using Bonferroni's inequality and taking into account relations~(\ref{v1})-(\ref{v2}), we have
\begin{eqnarray*}&&
P_n \geq E(W_1)+E(W_2)-E{W_1\choose 2}-E{W_2\choose 2} = \\& & \geq 4{n\choose 2}2^{-n} +16 {n\choose 4}\left(\frac{3}{8}\right)^n  -12{n\choose 2}^2 4^{-n} - 2^{16}n^5  \left(\frac{3}{16}\right)^n \\&& > 4{n\choose 2}2^{-n} +16 {n\choose 4}\left(\frac{3}{8}\right)^n  - 2^{17}n^5  \left(\frac{3}{16}\right)^n.\end{eqnarray*}
Conjecture~2 states that this lower bound is asymptotically tight. The history
of the problem of determining upper bound for $P_n$ started in 1967 when
in~\cite{2} Koml\'{o}s proved that $P_n = o(1)$. In 1995 in the work~\cite{7} Kahn, Koml\'{o}s and
Szemeredi proved that $P_n < (\alpha + o(1))^n$
for some $\alpha <1$  very closed to $1$.
Actually that work established many interesting ideas which also used later
in improvements of this bound. First such improvement was made in~\cite{3} by
Tao and Vu, where $\alpha$ was improved to $0.939$ and in the later work~\cite{4} to $0.75$. Their
improvement add additive combinatorics as ingredient in the proof. This last bound was improved by Bourgain,
Wood and Vu in~\cite{1} to 
\begin{equation}
\label{et1}
P_n \leq\left(\frac{1}{2}+o(1)\right)^{n/2}.
\end{equation}
In article~\cite{5} K.Tikhomirov proved tight logarithmic asymptotic
\begin{equation}
\label{b7}
P_n \leq \left(\frac{1}{2}+o(1)\right)^n.
\end{equation}
\bigskip

In this article we  prove  the following
\begin{Th}
\label{t3}
The following relation is valid
\begin{equation}
\label{21}
P_n \leq 4{n\choose 2} 2^{-n}+16{n\choose 4} \left(\frac{3}{8}\right)^n(1+o(1)).
\end{equation}
\end{Th}
Using the arguments in the proof of this theorem one can find the arbitrary fixed term expansion of $P_n$.

Proof of these Theorems allows to extend asymptotic  expansion of $P_n$ over $n$  with the arbitrary given precision.   
\bigskip

\newpage 
{\bf Proof}
\bigskip

At first we demonstrate rather short proof of the bound~(\ref{b7}), using results of previous work~\cite{1}. 

We divide linear subspaces  ${\cal R}$ of $ R^n$ into three   families  
\begin{eqnarray*}&&
{\cal R}_1 =\left\{  V\subset R^n: P(\{\pm 1\}^n \in V^\perp )> \frac{100}{\sqrt{n}}\right\}; \\
&&  {\cal R}_2=\left\{ V\subset R^n:\ (0.51)^{n/2}<P(\{\pm 1\}^n \in V^\perp)\leq \frac{100}{\sqrt{n}}\right\} ;\\&& {\cal R}_3=\left\{ V\subset R^n:\ P(\{\pm 1\}^n \in V^\perp)\leq (0.51)^{n/2}\right\},
\end{eqnarray*}
where $P(a)=2^{-n},\ a\in\{\pm 1\}^n$. 

Upper bound for the probability that the following statement is true $$\bigcup_{V\in {\cal R}_1}\{A_{n,n}\subset V^\perp\} \bigcup_{V\in {\cal R}_1}\{A^T_{n,n}\subset V^\perp\}$$ is
stated in~(\ref{bb11}). First inequality in~(\ref{bb11}) one can find in~\cite{7}. Note that if $v\in V\in {\cal R}_1$,  then  $|\hbox{Supp}_n (v)|<\frac{n}{2\cdot 10^4 \pi}< n-6\biggl[\frac{n}{\log_2 (n)}\biggr]$ (we demonstrate the proof  below for completeness see~(\ref{ss1})). Upper bound for the probability of event  $\bigcup_{V\in {\cal R}_2}\left\{A_{n,n}\subset V^\perp\right\}\bigcup_{V\in {\cal R}_2}\left\{A^T_{n,n}\subset V^\perp\right\} $  is stated in~(\ref{bb2})
(see proposition~5.4~\cite{1}).
Here event $\{A_{n,n}\in V^\perp\}$ means that $\{a_i,\ i\in [n]\}\subset V^\perp.$
   
 For $||x||=1,\ d=\min_{i:x_i\neq 0} |x_i |,\ \tilde{x}_i =\frac{x_i}{d},\ |\hbox{Supp}_n (x)|=k$ and $P(b_i=1)=P(b_i=-1)=\frac{1}{2}$ are independent variables, we have $(\epsilon =\frac{d}{2})$
\begin{eqnarray*}
&& P\left(|(x,b) -\lambda | <\epsilon \right) = P\left((x,b) \in \left(\lambda -\epsilon,\ \lambda+\epsilon\right) \right)\\
&& =P\left(( \tilde{x} ,b)\in \left(\frac{\lambda}{d}-\frac{\epsilon}{d},\frac{\lambda}{d}+\frac{\epsilon}{d}\right)\right) =P\left((\tilde{x},b) \in\left(\frac{\lambda}{d} -\frac{1}{2},\frac{\lambda}{d}+\frac{1}{2}\right)\right) \leq\frac{{n\choose[k/2]}}{2^k}.
\end{eqnarray*}
Last inequality is Erdos-Littlewood-Offord inequality~\cite{oo1}, which state that when $|\tilde{x_i}|\geq 1$ and $ |A|\leq 1,$ then
$$
P\left((b,\tilde{x}) \in A \right) \leq\frac{{n\choose[k/2]}}{2^k}.
$$
Hence  for sufficiently small $\omega >0$ we have
\begin{equation}
\label{r1}
\sup_{\lambda\in R}P(|(x,b) -\lambda| <\omega) \leq\frac{{k\choose[k/2]}}{2^k}.
\end{equation}
It follows
\begin{equation}
\label{k3}
P\left(|(x, b)| <\omega \right) \leq\frac{{k\choose[k/2]}}{2^k}.
\end{equation}
  When $P(\{\pm 1\}^n \in V^\perp)>\frac{100}{\sqrt{n}}$, then  it follows inequality $\min_{v\in V}P((a,v)=0)\geq \frac{100}{\sqrt{n}}$.
  
   Due to~(\ref{k3}) for given $|\hbox{Supp}_n (v)|=2k$ using relations

\begin{eqnarray}&&
\label{k9}
\frac{4^k}{\sqrt{\pi \left(k+\frac{1}{2}\right)}}\leq {2k\choose k}\leq \frac{4^k}{\sqrt{\pi k}},\\
&& \label{k10}   {2k+1\choose k}= {2k\choose k}\left(1+\frac{k}{k+1}\right) 
\end{eqnarray}
   
    we have  
  $$
  \frac{1}{\sqrt{2\pi k}}\geq \frac{{2k\choose k}}{2^{2k}}\geq \frac{100}{\sqrt{n}}.
  $$
  or
  \begin{equation}
  \label{ss1}
  k\leq \frac{n}{ 2\cdot 10^4 \pi}.\end{equation}
  Hence if $V\in {\cal R}_1$, then for arbitrary $v\in V$ we have $|\hbox{Supp}_n (v)| \in \left[ 2, n-6\left[\frac{n}{\log_2 (n)}\right]\right]$.

Condition
   $|\hbox{Supp}_n (v)|\leq n-6[\frac{n}{\log_2 (n)}]$ for some $v\in V\in {\cal R}_1$ or $|\hbox{Supp}_n (\bar{v})|\leq n-6[\frac{n}{\log_2 (n)}]$ for some $\bar{v}\in V\in {\cal R}_1$ is sufficient for the inequality 
  \begin{eqnarray}&&
  \label{bb11}
  P\left(\bigcup_{V\in {\cal R}_1}\{A_{n,n}\in V^\perp\}\bigcup_{V\in {\cal R}_1}\{A_{n,n}^T\in V^\perp\}\right) \leq 4{n\choose 2}2^{-n}+16{n\choose 4}\left(\frac{3}{8}\right)^n\\&& +2\sum_{k=5}^{n-6[\frac{n}{\log_2 (n)}]} {n-1\choose k-2}{n\choose k} p_k^{n-k+1} =4{n\choose 2}2^{-n}+ 16{n\choose 4}\left(\frac{3}{8}\right)^n(1+o(1)).\nonumber
   \end{eqnarray}   to  be true. Here $p_k=\frac{{k\choose [k/2]}}{2^k}.$ 
   The expression in the middle of chain of inequalities~(\ref{bb11})  arise from the following consideration: we can choose pairs of columns or rows 2${n\choose 2}$ possible ways and
   the probability of linear dependence of rows or columns is $2^{1-n}$, the same consideration for fours of columns or rows leads to second term. The third term - sum arise from the consideration that  we can choose submatrix $A_{k,k-1}$ of rank $k-1$ from matrix $A_{n,n}$ by ${n\choose k}{n-1\choose k-2}$ ways (second binomial coefficient arise from the fact that we can fix first column of $A_{n,n}$). Then the probability of any choice of other $n-k+1$ columns of matrix $A_{k,n}$ has probability less that $p_k^{n-k+1}$. We make the same consideration for $A_{n,n}^T$.    At last we take the sum over all $k\in \left[ 2, n-6\left[\frac{n}{\log_2 (n)}\right]\right]$.
         
  When $V\in{\cal R}_2$ we have $(0.51)^{n/2} <P(\{\pm 1\}^n \in V^\perp )<\frac{100}{\sqrt{n}}$ and we use bound~(Proposition~5.4,~\cite{1}):
  \begin{equation}
  \label{bb2}
  P\left(\bigcup_{V\in {\cal R}_2}\{A_{n,n} \in V^\perp\}\right) \leq (o(1))^n.
  \end{equation}
 Denote the set of linear spaces $W=\{V \in {\cal R}_3 ,\ \hbox{dim}(V)\leq n-3\}$. 
 
 Then for large $n$
\begin{eqnarray}&&
\label{tr3}
P\left(\bigcup_{V\in W} \{A_{n,n}\in V^\perp\}\right) \leq \sum_{i=1}^{n-3}{n\choose i}(0.51)^{n(n-i)/2}\\ &&= ((0.51)^{n/2} +1)^n -1-n(0.51)^{n/2}-{n\choose 2}(0.51)^{n}<n^3(0,51)^{3n/2}.\nonumber
\end{eqnarray} 
It follows from~(\ref{bb11}),~(\ref{bb2}),~(\ref{tr3}) that 
\begin{eqnarray}&&
\label{h1}
 P\Biggl( \bigcup_{V\in W} \{A_{n,n}\in V^\perp\}\bigcup_{V\in {\cal R}_1}\{A_{n,n}\in V^\perp\}\bigcup_{V\in {\cal R}_1}\{A^T_{n,n}\in V^\perp\}\bigcup_{V\in {\cal R}_2}\{A_{n,n}\in V^\perp\}\Biggr)\\&&  \nonumber   \leq 4{n\choose 2} 2^{-n} +16{n\choose 4}\left(\frac{3}{8}\right)^n(1+o(1)).
\end{eqnarray}
If $P(\{\pm 1\}^n \in  V^\perp)>\frac{100}{\sqrt{n}}$, then for all $v\in V$ we have $ |\hbox{Supp}_n (v) |\leq  n-6[\frac{n}{\log_2 (n)}]$. 

Hence we need to consider the case when $|\hbox{Supp}_n (v)|>n-6\left[\frac{n}{\log_2 (n)}\right]$ for all $v(A_{n,n})=(v_1,\ldots ,v_n)\in V\in \cal{R}_3$.

Denote $r=n-m=[0.7n]$ and $|\hbox{Supp}_{n-m} (v(A_{n,n}))|=p,\ \max_{|v_{i}|\neq 0,i\leq n-m}i=q$;
\begin{eqnarray}\nonumber &&\bar{v}_{b,r}(A_{n,n})=\frac{1}{\sqrt{\sum_{i=1}^{n-q-1} v_i^2 (A_{n,n})+\bar{v}_{n-q}^2 (A_{n,n})}}(v_1 (A_{n,n}), \ldots , v_{n-q-1}(A_{n,n}), \bar{v}_{n-q}(A_{n,n}), 0,\ldots , 0);\\&& b=(b_{r+1},\ \ldots ,b_n)\in \cal{B}=\{\pm 1\}^{m},\ \bar{v}_{n-q}(A_{n,n}) = v_{n-q}(A_{n,n})+\sum_{i=r+1}^n b_i v_i (A_{n,n}). \label{rr1}\end{eqnarray} For every $v(A_{n,n}),\ b\in{\cal B}=\{\pm 1\}^{m}$ we have
\begin{eqnarray}&& \label{ttr}
P((a,\bar{v}_{b,r}(A_{n,n}))=0)\leq P((a,v(A_{n,n}))=0) \leq \sum_{b\in\cal{B}}P((a,\bar{v}_{b,r}(A_{n,n}))=0)\\ &&\leq 2^{m}\max_{b\in\cal{B}}P((a,\bar{v}_{b,r}(A_{n,n})=0)\nonumber
\end{eqnarray}
and 
$$
\{a:\ P((a,v(A_{n,n}))=0\}=\bigcup_{b\in\cal{B}} \{a:\ P((a,\bar{v}_{b,r}(A_{n,n}))=0)\} .
$$

When $\dim{V}=n-1, v\in {\cal R}_3$  and $
|\hbox{Supp}_n (v)| >n-6\left[\frac{n}{\log_2 (n)}\right]$, then for $r=n-m=[0.7n]$ we have
\begin{eqnarray}&&
  P\left(\bigcup_{V\in {\cal R}_3}\{A_{n,n}\in v^\perp\}\bigcup_{v\in {\cal R}_3}\{A_{n,n}^T\in v^\perp\}\right) \label{jj1}
   \leq P\left(\bigcup_{v\in {\cal R}_3}\{A_{n,n}\in v^\perp\}\bigcup_{v\in {\cal R}_3}\{A_{n,n}^T\in v^\perp\}\right)\\&& \nonumber
   \leq 2 \sum_{d=n-6\left[\frac{n}{\log_2 (n)}\right]}^{n}{n-1\choose d-2}{n\choose d}{n\choose [0.7n]}((0.51)^n n2^{[0.3n]} )^{[0.3n]} < 2^{3n} ((0.51)^n n2^{[0.3n]} )^{[0.3n]} =(o(1))^n. \nonumber\end{eqnarray}
The case $d\leq n-6\left[\frac{n}{\log_2 (n)}\right]$ is completely considered in~(\ref{bb11}). Considering other case $d> n-6\left[\frac{n}{\log_2 (n)}\right]$ we have that for each choice of $d\times (d-1)$
submatrix of matrix $A_{n,n}$ generate vector $v\in {\cal R}_3$ of with $\hbox{Supp}_n(v)$ nonzero coordinates, which according to~(\ref{rr1}) generate the set of vectors $\bar{v}_{b,r}$ of length $r$  with the property~(\ref{rr1}) and hence the set of matrices s.t. $A_{n,n}\in v^\perp$ covered by the set of matrices, generated by the set of relations $A_{n,n}\in \bar{v}^\perp_{b,r}$. Hence each choice of submatrix $C_{r,r-1}$ of size $r\times (r-1)$ from matrix $A_{n,n}$ allows to generate the set of possible $r$ column vectors $\bar{c}_i =(c_{1,i},\ldots , c_{r,i}),\ i\in [r,n]$ in matrix $C_{r,n}$,  each with probability at most $2^{m}\max_{b\in\cal{B}}P((a,\bar{v}_{b,r})=0)\leq 2^m (0.51)^n$. 

It is left to consider the case 
\begin{eqnarray*}&&
Q=\left\{V=(v_1,v_2),\ \hbox{dim}(V)=n-2,\  |\hbox{Supp}_n (v_1)|,\ |\hbox{Supp}_n (v_2)| >n-6\left[\frac{n}{\log_2 (n)}\right]\right\}\subset{\cal R}_3. 
\end{eqnarray*}
Assume that $r=n-m=[0.7n]$ and $|\hbox{Supp}_{n-m} (v_j)|=p_j ,\ \max_{|v_{j,i}|\neq 0,\ i\leq r}i=q_j,\ j=1,2$.    

Denote
\begin{eqnarray*}&&\bar{v}_{1,b,r}(A_{n,n})=\frac{1}{\sqrt{\sum_{i=1}^{q_1-1} v_{1,i}^2(A_{n,n}) +\bar{v}_{1,q_1}^2 (A_{n,n})}}(v_{1,1} (A_{n,n}), \ldots , v_{1,q_1-1}(A_{n,n}), \bar{v}_{1,q_1}(A_{n,n}), 0,\ldots , 0);\\
&&\bar{v}_{2,b,r}(A_{n,n})=\frac{1}{\sqrt{\sum_{i=1}^{q_2-1} v_{1,i}^2(A_{n,n}) +\bar{v}_{2,q_2}^2 (A_{n,n})}}(v_{2,q_2} (A_{n,n}), \ldots , v_{2,q_2-1}(A_{n,n}), \bar{v}_{2,q_2}(A_{n,n}), 0,\ldots , 0);\\
&&b=(b_{r+1},\ \ldots ,b_n)\in \cal{B}=\{\pm 1\}^{m},\ \bar{v}_{1,q_1}(A_{n,n}) = v_{1,q_1}(A_{n,n})+\sum_{i=r+1}^n b_i v_{1,i} (A_{n,n});\\
&&\bar{v}_{2,q_2}(A_{n,n}) = \bar{v}_{2,q_2}(A_{n,n})+\sum_{i=r+1}^n b_{i} \bar{v}_{2,i} (A_{n,n}).
\end{eqnarray*} For every $(v_1(A_{n,n}),v_2(A_{n,n}))\in {\cal R}_3,\ b\in\cal{B}=\{\pm 1\}^{m}$ we have
\begin{eqnarray}&& \label{ttr}
P((a,\bar{v}_{1,b,r}(A_{n,n}))=0)\leq P((a,v(A_{n,n}))=0) \leq \sum_{b\in\cal{B}}P((a,\bar{v}_{1,b,r}(A_{n,n}))=0)\\ &&\leq 2^{m}\max_{b\in\cal{B}}P((a,\bar{v}_{1,b,r}(A_{n,n}))=0);\nonumber\\
&&P((a,\bar{v}_{2,b,r}(A_{n,n}))=0)\leq P((a,v(A_{n,n}))=0) \leq \sum_{b\in\cal{B}}P((a,\bar{v}_{2,b,r}(A_{n,n}))=0)\nonumber \\ &&\leq 2^{m}\max_{b\in\cal{B}}P((a,\bar{v}_{2,b,r}(A_{n,n}))=0)\nonumber
\end{eqnarray}
and 
\begin{eqnarray*}
&&
\{a:\ P((a,v_1(A_{n,n}))=0)\}=\bigcup_{b\in\cal{B}} \{a:\ P((a,\bar{v}_{1,b,r}(A_{n,n}))=0)\}; \\
&& \{a:\ P((a,v_2(A_{n,n}))=0)\}=\bigcup_{b\in\cal{B}} \{a:\ P((a,\bar{v}_{2,b,r}(A_{n,n}))=0)\}. 
\end{eqnarray*} 
 For $r=n-m =[0.7n]$ we have with the proof simular to~(\ref{jj1}) the following bound
   \begin{eqnarray}&&
  \label{f11}
  P\left(\bigcup_{V=\{v_1,v_2\}\in {\cal R}_3}\{A_{n,n}\in V^\perp\}\bigcup_{V=\{v_1,v_2\}\in {\cal R}_3}\{A_{n,n}^T\in V^\perp\}\right)\nonumber
  \\&& \nonumber
   \leq 2 \sum_{d=n-6\left[\frac{n}{\log_2 (n)}\right]}^{n}{n-1\choose d-2}{n\choose d}{n\choose [0.3n]}((0.51)^n n2^{2[0.3n]} )^{2[0.3n]} < 2^{3n} ((0.51)^n n2^{2[0.3n]} )^{[0.3n]} =(o(1))^n. \nonumber\end{eqnarray}
Proof of Theorem~2 is completed.

\newpage

\begin{center}
{\bf Acknowledgment}
\end{center}

Research was started under support by the Sao Paulo Research Foundation (FAPESP), projects no. 2012/13341-8, 2014/23368-6 and NUMEC/USP project 2013/07699-0 and visiting position in Unifesp
by Extrato de contrato N114/2019, N20/2021, processo N23089.101560/2018-91.  We would like to express our gratitude to  Unifesp  and also to USP,  where first author started this work and specially to Prof.Y.Kohayakawa.    

Research supported by Extrato de contrato N114/2019, N20/2021, Processo N23089.101560/2018-91.  The article was
written while the author was visiting Unifesp, Brazil. He would like to
thank the Department of Mathematics and also colleagues for their kind atmosphere to execute  this work.

\end{document}